\documentclass[12pt]{article}

\usepackage{amssymb}
\usepackage{amsmath}
\usepackage{amsfonts}
\usepackage{amsthm}
\usepackage{epsfig}

\title{Overcrowding estimates for zeroes \\  of Planar and Hyperbolic
  \\ Gaussian analytic functions} 
\author{Manjunath Krishnapur\footnotemark[1]}
  \footnotetext[1]{Research supported by NSF 
  grant \#DMS-0104073 and NSF-FRG grant \#DMS-0244479.}

\theoremstyle{theorem}
    \newtheorem{theorem}{Theorem}
    \newtheorem{lemma}[theorem]{Lemma}

\theoremstyle{definition} 

    \newtheorem{remark}[theorem]{Remark}

\theoremstyle{remark}

\def\tends{\rightarrow}

\def\P{{\bf P}}

\def\eqd{\stackrel{d}{=}}

\def\Gam{\Gamma}

\def\alp{\alpha}
\def\bet{\beta}
\def\gam{\gamma}
\def\del{\delta}

\def\eps{\epsilon}
\def\lam{\lambda}

\def\zet{\zeta}

\def\Mid{\left.\vphantom{\hbox{\Large (}}\right|}
\def\abs
\def\bar{\overline}
\def\summ{\sum\limits}
\def\prodd{\prod\limits}
\def\intt{\int\limits}
\def\l{\left}
\def\r{\right}

\def\hsp{\hspace}
\def\mb{\mbox}

\def\ind{{\mathbf 1}}

\def\d{\partial}

\def\<{\langle}
\def\>{\rangle}

\def\f{{\bf f}}
\def\g{{\bf g}}
\def\T{{\mathbb T}}

\begin{document}
\maketitle 
\begin{abstract} We consider the point process of zeroes of certain
  Gaussian analytic functions and find the asymptotics for the
  probability that there are more than $m$ points of the process in a
  fixed disk of radius $r$, as $m\tends \infty$. For the Planar
  Gaussian analytic function, $\sum_{n\ge0} \frac{a_n
  z^n}{\sqrt{n!}}$, we show that this probability is asymptotic to
  $e^{-\frac{1}{2}m^2\log(m)}$. For the Hyperbolic Gaussian analytic
  functions, $\sum_{n\ge 0}{-\rho \choose n}^{1/2} a_n z^n$,
  $\rho>0$, we show that this probability decays like $e^{-cm^2}$.  
 
 In the planar case, we also consider the problem posed by  Mikhail
 Sodin~\cite{sod2} on moderate and very large deviations  in a 
 disk of radius 
 $r$ as $r\tends \infty$. We partially solve the problem by showing that
 there is a qualitative change in the asymptotics of the probability
 as we move from the large deviation regime to the moderate.  
\end{abstract}
 
\section{Introduction} In this paper we consider the following
 Gaussian analytic functions (GAFs): 
\begin{itemize}
\item {\bf Planar GAF :} Often called the {\it Chaotic analytic function} in the Physics literature, this is the random analytic function
\[ \g(z)=\summ_{n=0}^{\infty} \frac{a_n z^n}{\sqrt{n!}} \]
where $a_n$ are i.i.d. standard Complex Normal random variables. This
defines an entire function (almost surely).
\newpage
\item {\bf Hyperbolic GAFs :} For each $\rho>0$ let
\[ \f_{\rho}(z) = \summ_{n=0}^{\infty}{-\rho \choose n}^{1/2} a_n z^n \]
where as before $a_n$ are i.i.d.  standard Complex Normals. Almost surely, $\f_{\rho}$ is an analytic function in the unit disk (and no more).  
\end{itemize}

These particular GAFs are of interest because the distributions of
their zero sets are invariant under isometries of the Euclidean plane
and isometries of the Hyperbolic plane respectively. In particular the
zero set of $\f_{\rho}$ has constant intensity $\frac{\rho}{\pi}$
(w.r.t. $\frac{d|z|^2}{(1-|z|^2)^2}$) and is the only zero set of a
GAF that is conformally invariant in the unit disk and has this
density. See Sodin and Tsirelson~\cite{sodtsi1} and Sodin~\cite{sod}
for proofs of these assertions. In the planar case too one can define
GAFs with invariant zero distribution of intensity $\rho$ for any
$\rho>0$, but these are just scaled versions of the zero set of $\g$
defined above. The zero set of $\g$ has intensity $\frac{1}{\pi}$
w.r.t. Lebesgue measure on the plane(and again is the only GAF zero
set with this intensity).  

 We denote the zero set by ${\mathcal Z}$. Let $n(r)$ denote the
 number of points of ${\mathcal Z}$ in the disk of radius $r$ around
 $0$ (The GAF will be clear from the context). We address the following
 two problems.
 \begin{enumerate}
 \item {\bf Overcrowding: } Yuval Peres asked the following question
 and conjectured that the probability decays as $e^{-cm^2\log(m)}$ in
 the planar case (personal communication). 
 
{\bf Question:} Fix $r>0$, ($r<1$ in the Hyperbolic case). Estimate
$\P\l[n(r)>m \r]$ as $m\tends \infty$. 
\begin{figure}
\includegraphics[height=2.8in]{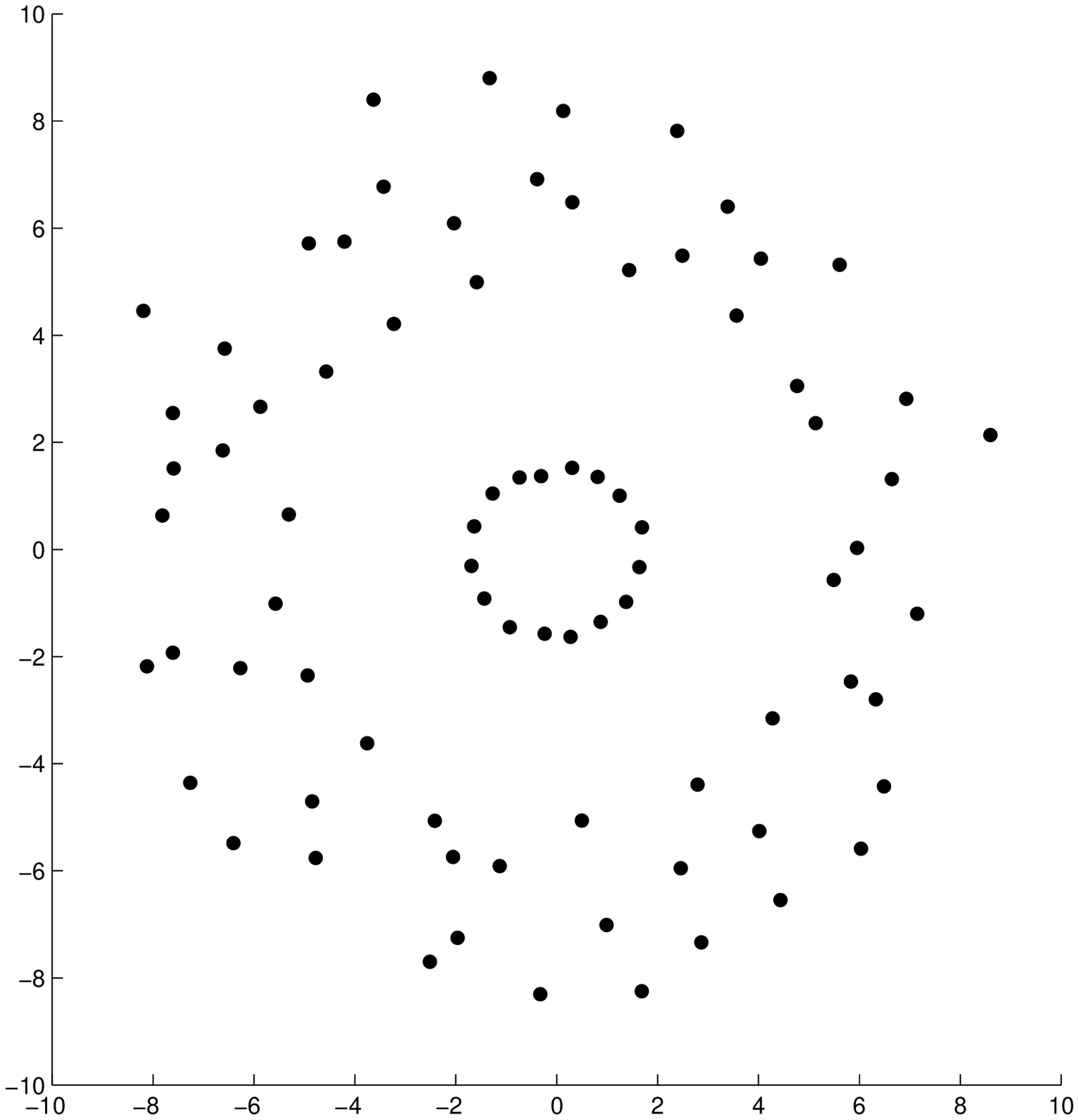}\hspace{.25in}
\includegraphics[height=2.8in]{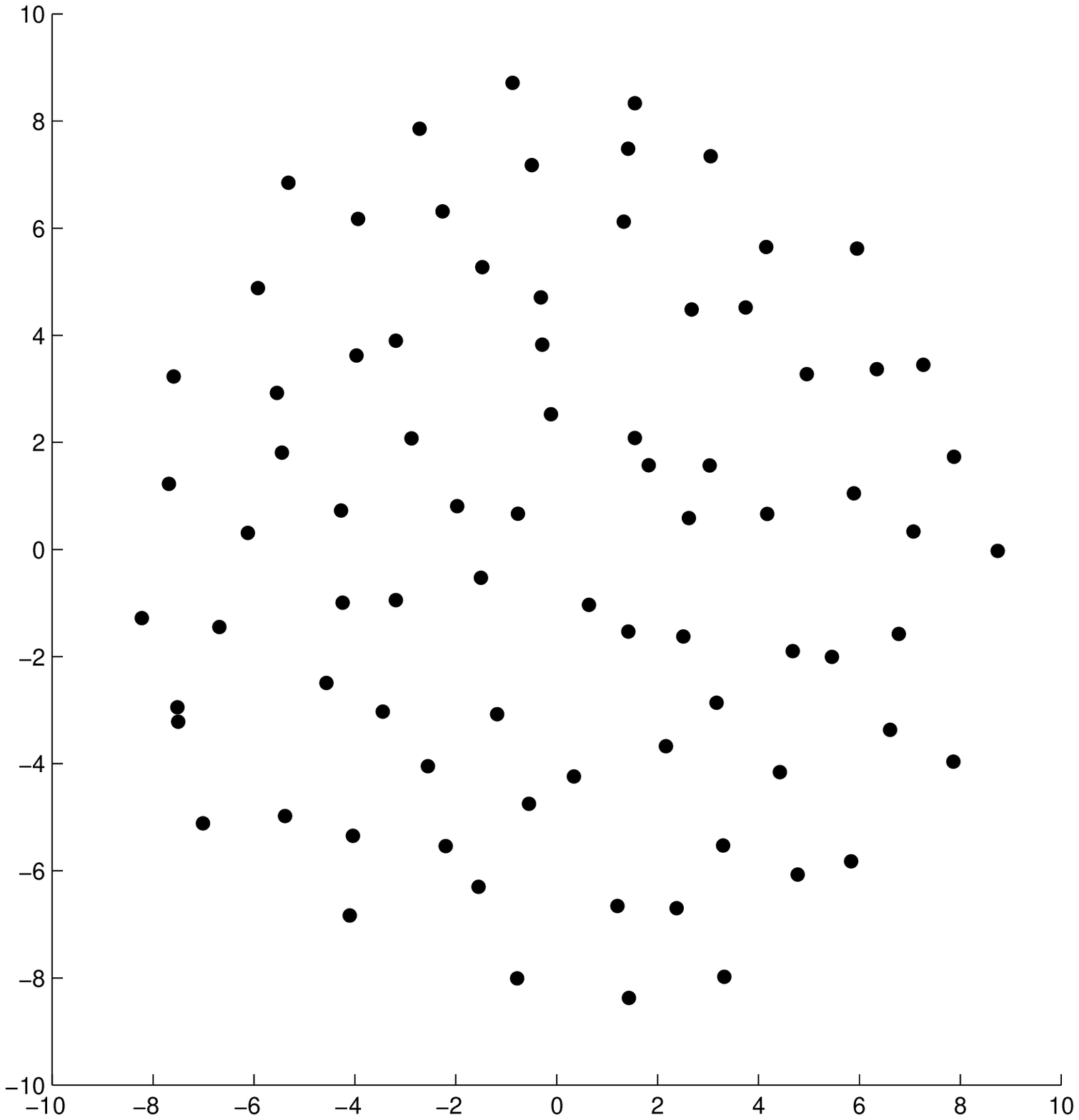}
  \centering
  \caption{Samples of the zero process of $\g$. Left: The zero process
    sampled under {\it certain sufficient
    conditions} (see the conditions in the lower bound of the proof of
    Theorem~\ref{thm:planar}. Take $\alp=0.5,r=2,m=16$) on the
    coefficients forcing $16$  zeroes in the disk of radius
    $2$. Right: The unconditioned zero process.} 
  \label{fig:unconditioned}
\end{figure}

One motivation for such a question is in
Figure~\ref{fig:unconditioned}. There one can see the distribution of the 
zero process under certain conditions on the coefficients that force 
large number of zeroes in the disk of radius $2$ (this is {\it not}
the zero set conditioned to have overcrowding - that seems harder to
simulate). The picture suggests that the distribution of the 
conditioned process may be worth studying on its own. A large
deviation estimate of the kind we derive will presumably be a
necessary step in such investigations.

The answer is different in the two settings. We prove-

\begin{theorem}\label{thm:planar} Consider the planar GAF $\g$. 
For any $\eps >0$, $\exists$ a constant $C_2$ (depending on $\eps,r$)
such that for every $m\ge 1$, 
\begin{equation*}
 e^{-\frac{1}{2} m^2\log(m)+O(m^2)} \le \P[n(r)\ge m] \le C_2
 e^{-(\frac{1}{2}-\eps)m^2\log(m)}. 
\end{equation*} 
In particular, $\P[n(r)\ge m]=e^{-\frac{1}{2} m^2\log(m)(1+o(1))}$.
\end{theorem}

\begin{theorem}\label{thm:hyperbolic} Fix $\rho>0$ and consider the
  GAF  $\f_{\rho}$. For any fixed $r<1$, there are constants
  $\bet,C_1,C_2$(depending 
  on $\rho$ and $r$) such that for every $m\ge 1$, 
\begin{equation*} C_1(r)e^{- \frac{m^2}{|\log(r)|}} \le
  \P[n(r)\ge m] \le C_2(r) e^{-\bet(r) m^2}.
\end{equation*}
\end{theorem}

\item {\bf Moderate, Large and Very Large Deviations:}
Inspired by the results obtained by Jancovici, Lebowitz and 
Manificat~\cite{janlebmag} for Coulomb gases in the plane (eg.,
Ginibre ensemble), M.Sodin~\cite{sod2} has conjectured the following. 

{\bf Conjecture: } Let $n(r)$ be the number of zeroes of the planar
GAF $\g$ in the disk $D(0,r)$. Then, as $r\tends \infty$ 
\begin{equation}\label{eq:conjecture}
\frac{\log \log \l(\frac{1}{\P[|n(r)-r^2|>r^{\alp}]}\r)}{\log r} \tends \l\{ \begin{array}{cc}
         2\alp -1, & \frac{1}{2}\le \alp \le 1; \\
         3\alp -2, & 1\le \alp \le 2; \\
         2\alp,    & 2\le \alp. \end{array} \r. 
\end{equation}
The idea here is that the deviation probabilities undergo a
qualitative change in behaviour when the deviation under consideration
becomes comparable to the perimeter ($\alp=1$) or to the area
($\alp=2$) of the domain.

Sodin and Tsirelson~\cite{sodtsi3} had already settled the case
$\alp=2$ by showing that for any $\del>0$, $\exists
c_1(\del),c_2(\del)$ such that 
\begin{equation*}
  e^{-c_1(\del)r^4} \le \P[|n(r)-r^2|>\del r^2] \le e^{-c_2(\del)r^4}.
\end{equation*}

Here we consider $\P[n(r)-r^2>r^{\alp}]$ and prove that a phase
transition in the exponent occurs at  $\alp=2$. More precisely we
prove that the conjecture holds for $\alp>2$ and show the lower bound
for $1<\alp<2$. 

\begin{theorem}\label{thm:verylarge} Fix $\alp>2$. Then 
\[ \P\l[n(r)\ge r^2+\gam r^{\alp} \r] = e^{-\l(\frac{\alp}{2} -
  1\r)\gam^2 r^{2\alp}\log r (1+o(1))}. \] 
\end{theorem}

 \begin{theorem}\label{thm:moderate} Fix $1<\alp < 2$. Then  for
   any $\gam>0$, 
\[ \P\l[n(r)\ge r^2+\gam r^{\alp} \r]\ge e^{-\gam^3 r^{
    3\alp-2}(1+o(1))}. \] 
\end{theorem}


\begin{remark}  Nazarov, Sodin and Volberg have recently proved all
  parts of the conjecture (personal communication).
\end{remark}
\end{enumerate}

We prove Theorem~\ref{thm:planar} in Section~\ref{sec:planar},
Theorem~\ref{thm:hyperbolic} in Section~\ref{sec:hyperbolic}, and
Theorem~\ref{thm:verylarge} and Theorem~\ref{thm:moderate} in Section~\ref{sec:moderate}.

\section{Overcrowding - The Planar case}\label{sec:planar}
In this section we prove Theorem~\ref{thm:planar}. Before that we
explain  why one expects the constant $\frac{1}{2}$ in the exponent in
Theorem~\ref{thm:planar}, by analogy with the Ginibre ensemble. 
 
\subsection{Ginibre Ensemble} 
The Ginibre ensemble is the determinantal point process (see
\cite{sos} or \cite{hkpv} for definitions) in the plane with kernel
\begin{equation}\label{eq:ginkernel}
 K(z,w)=\frac{1}{\pi}e^{-\frac{1}{2}|z|^2-\frac{1}{2}|w|^2+z{\bar w}}. 
\end{equation}
This process is of interest because it is the limit in distribution,
 as $n\tends\infty$, of the  
 point process of eigenvalues of an $n\times n$ matrix with
 i.i.d. standard Complex Normal entries~\cite{gin}. 

The Ginibre ensemble has many similarities to the zero set of $\g$. In
particular, the Ginibre ensemble is invariant in distribution under
Euclidean motions, 
has constant intensity $\frac{1}{\pi}$ in the plane and has the same
negative correlations as ${\mathcal Z}_{\g}$ at short
distances. Therefore there are other 
similarities too, for instance, see~\cite{denhan}. There are also
differences between the two point processes. For
instance, the Ginibre ensemble has all correlations negative, whereas
for the zero set of $\g$, long-range two-point correlations are
positive. However, in our problem, since we are considering a fixed
disk and looking at the event of having an excess of zeroes in it, it
seems reasonable to expect the same behaviour for both these point
processes, since it is the short range interaction that is
relevant. In case of the Ginibre ensemble, the overcrowding problem is
easy to solve. 

\begin{theorem}\label{thm:ginibre}
Let $n_G(r)$ be the number of points of the Ginibre ensemble in the disk of radius $r$ around $0$ (by translation invariance, the same is true for any disk of radius $r$). Then for a fixed $r>0$,
\[ \P\l[n_G(r)\ge m\r] = e^{-\frac{1}{2}m^2\log(m) (1+o(1))}. \]
\end{theorem}
\begin{proof} By Kostlan~\cite{kostlan}, the set of absolute values of
  the points of the Ginibre ensemble has the same distribution as the
  set $\{R_1,R_2,\ldots \}$, where $R_n$ are independent, and $R_n^2$
  has Gamma($n,1$) distribution for every $n$. Hence $R_n^2\eqd
  \xi_1+\ldots +\xi_n$, where $\xi_k$ are i.i.d.  Exponential random
  variables with mean $1$, and it follows that 
\[ \P\l[R_n^2<r^2 \r]\ge \prodd_{k=1}^n \P\l[\xi_k<\frac{r^2}{n}\r]
  \ge \l(\frac{r^2}{2n}\r)^n, \] as long as $n\ge r^2$, because
  $\P\l[\xi_1<x\r]\ge \frac{x}{2}$ for $x<1$. Therefore we get
\begin{eqnarray} 
\P\l[n_G(r)\ge m\r] &\ge& \prodd_{n=1}^{m} \P\l[R_n^2<r^2 \r] \\
 &\ge& \prodd_{n=1}^{m} \l(\frac{r^2}{2n}\r)^n \\
 &=& \l(\frac{r^2}{2}\r)^{\frac{m(m+1)}{2}}e^{-{\summ_{n=1}^m
 n\log(n)}}.
\label{eq:lbfornr}
\end{eqnarray}
Here and elsewhere we shall encounter the term $\summ_{n=1}^m
n\log(n)$. We compute its asymptotics now. 
\[ n\log(n) \le x\log(x) \le (n+1)\log(n+1) \hsp{2cm}\mb{ for }n\le x\le n+1 \]
Integrate from $1$ to $m+1$ and note that 
\[ \intt_1^a x\log(x) dx = \frac{1}{2}a^2 \log(a) -
\frac{a^2}{4}+\frac{1}{4}, \] 
to get
\begin{equation}\label{eq:sumnlogn}
\summ_{n=1}^m n\log(n)\le  \frac{1}{2}(m+1)^2 \log(m+1) - \frac{(m+1)^2}{4}+ \frac{1}{4} \le \summ_{n=1}^{m+1} n\log(n). 
\end{equation}
Thus (\ref{eq:lbfornr}) gives
\begin{eqnarray*}
\P\l[n_G(r)\ge m\r] &\ge& e^{-\frac{1}{2}(m+1)^2\log(m+1)+\frac{(m+1)^2}{4}-\frac{1}{4}+\frac{m(m+1)}{2}\log(r^2/2)}\\
 &=& e^{-\frac{1}{2}m^2\log(m)+O(m^2)}.
\end{eqnarray*}
To prove the inequality in the other direction, note that
\begin{eqnarray*}
\P\l[n_G(r)\ge m\r] &\le& \P\l[\summ_{n=1}^{m^2}\ind(R_n^2<r^2)\ge m \r]+\summ_{n=m^2+1}^{\infty} \P\l[R_n^2<r^2\r] \\
&\le& {m^2 \choose m}\prodd_{n=1}^{m} \P\l[R_n^2<r^2 \r] +
\summ_{n>m^2}e^{-n\log(n)(1+o(1))}.
\end{eqnarray*}
In the second line, for the first summand we used the fact that $R_n^2$ are
stochastically increasing and for the second term we used the well-known
fact $\P\l[R_n^2<r^2\r]=\P\l[\mb{Pois}(r^2)\ge n\r]$ and then the
usual bound on the tail of a Poisson random variable, namely $\P\l[\mb{Poisson}(\theta) \ge a \r]\le e^{-a\log(a/\theta)+a-\theta}$. 

Using the same idea to bound $\P\l[R_n^2<r^2\r]$ in the first summand,
we obtain
\begin{eqnarray*}
\P\l[n_G(r)\ge m\r] &\le& {m^2 \choose m}\prodd_{n=1}^{m} e^{-{n\log(n/r^2)-r^2+n}}+ e^{-m^2\log(m^2)(1+o(1))} \\
&\le& {m^2 \choose m}e^{\frac{m(m+1)}{2}(1+\log(r^2))-mr^2-\summ_{n=1}^m n\log(n)}+e^{-m^2\log(m^2)(1+o(1))} \\
&=& e^{-\frac{1}{2}m^2\log(m) (1+o(1))} \hsp{2cm} \l(\mb{using (\ref{eq:sumnlogn}) again} \r).
\end{eqnarray*}
In the last line we used ${m^2\choose m}<m^{2m}$. This completes the proof. 
\end{proof}

\subsection{Proof of Theorem~\ref{thm:planar}}
 Our method of proof is largely based on  that of Sodin and
 Tsirelson~\cite{sodtsi3}. 
(They estimate the ``hole probability'', $\P\l[ n(r)=0\r]$ as $r\tends
 \infty$.) 

\begin{proof}[Proof of Theorem~\ref{thm:planar}]

{\bf Lower Bound } Suppose the $m^{\mb{th}}$ term dominates the sum of all the  other terms on $\d D(0;r)$, i.e., suppose
\begin{equation}\label{eq:dominate} 
\Mid\frac{a_m z^m}{\sqrt{m!}}\Mid \ge \Mid \summ_{n\neq m} \frac{a_n z^n}{\sqrt{n!}}\Mid \hspace{2cm} \mb{ whenever }|z|=r.
\end{equation}
Then, by Rouche's theorem $\g(z)$ and $\frac{a_m z^m}{\sqrt{m!}}$ have the same number of zeroes in $D(0;r)$. Hence $n(r)=m$. 
Now we want to find a lower bound for the probability of the event in (\ref{eq:dominate}). Note that the left side of (\ref{eq:dominate}) is identically equal to $\frac{|a_m|r^m}{\sqrt{m!}}$. 

Now suppose the following happen-
\begin{enumerate}\label{enum:conditions}
\item $|a_n| \le n$ $\forall n\ge m+1$.
\item $|a_m| \ge (\alp+1)m$ where $\alp$ will be chosen shortly.
\item $|a_n|\frac{r^{n}}{\sqrt{n!}} < \frac{r^m}{\sqrt{m!}}$ for every $0\le n\le m-1$.
\end{enumerate}

Then the right hand side of (\ref{eq:dominate}) is bounded by
\begin{eqnarray*}
\mb{RHS of } (\ref{eq:dominate}) &\le& \summ_{n=0}^{m-1}|a_n|\frac{r^{n}}{\sqrt{n!}}  + \summ_{n=m+1}^{\infty} \frac{|a_n|r^n}{\sqrt{n!}} \\
&\le& \summ_{n=0}^{m-1}\frac{r^{m}}{\sqrt{m!}} + \summ_{n=m+1}^{\infty} \frac{nr^n}{ \sqrt{n!}} \\
&\le& m\frac{r^{m}}{\sqrt{m!}} + C \frac{m r^m}{\sqrt{m!}} \\
&=& (C+1)m\frac{r^{m}}{\sqrt{m!}} \\
&\le& |a_m|\frac{r^m}{\sqrt{m!}}
\end{eqnarray*}
if $\alp=C$. Thus if the above three events occur with $\alp=C$, then the $m^{\mb{th}}$ term dominates the sum of all the other terms on $\d D(0;r)$. Also these events have probabilities as follows.
\begin{enumerate}
\item $\P[|a_n| \le n$ $\forall n\ge m+1]\ge 1- \summ_{n=m+1}^{\infty} e^{-n^2}\ge 1-C'e^{-m^2}$. 
\item $\P[|a_m|\ge (C+1)m] = e^{-(C+1)^2m^2}$.
\item The third event has probability as follows. Recall again that
  $\P\l[\xi<x\r]\ge \frac{x}{2}$ if $x<1$ and $\xi$ is Exponential with
  mean $1$. We apply this below with $x=\l(\frac{r^{m-n}\sqrt{n!}}{
    \sqrt{m!}}\r)^2$. This is clearly less than $1$ if $n\ge
  r^2$. Therefore if $m$ is sufficiently large it is easy to see that
  for all $0\le n\le m-1$, the same is valid. Thus
\begin{eqnarray*}
\P\l[|a_n|\le \frac{r^{m-n}\sqrt{n!}}{\sqrt{m!}} \hsp{2mm}\forall n\le m-1 \r] &=& \prodd_{n=0}^{m-1} \P\l[|a_n|\le \frac{r^{m-n}\sqrt{n!}}{\sqrt{m!}} \r] \\
&\ge& \prodd_{n=0}^{m-1} \frac{r^{2m-2n}n!}{2 m!} \\
&=& r^{m(m+1)}e^{\frac{1}{2}m^2\log(m)+O(m^2)} 2^{-m} e^{-m^2\log(m)+O(m^2)} \\
&=& e^{-\frac{1}{2}m^2\log(m)+O(m^2)}.
\end{eqnarray*}
\end{enumerate}
Since these three events are independent, we get the lower bound in the theorem.

{\bf Upper Bound } By Jensen's formula, for any $R>r$ we have
\begin{equation}\label{eq:jensens}
n(r)\log\l(\frac{R}{r}\r) \le \intt_r^{R} \frac{n(u)}{u} du = \intt_0^{2\pi}\log|\g(Re^{i\theta })|\frac{d\theta}{2\pi} - \intt_0^{2\pi}\log|\g(re^{i\theta })|\frac{d\theta}{2\pi}.
\end{equation}

Let $R=R_m=\sqrt{m}$. Sodin and Tsirelson~\cite{sodtsi3} show that 
\begin{equation}\label{eq:ubdformax}
\P\l[\log M(t)\ge\l(\frac{1}{2}+\eps\r)t^2\r]\le e^{-e^{\eps t^2}}
\end{equation}
where $M(t)=\max\{|\g(z)|: |z|\le t \}$. 

Now suppose $n(r)\ge m$ and $\log M(R_m)\le \l(\frac{1}{2}+\eps \r)m$ for some $\eps>0$. Then by (\ref{eq:jensens}) we have
\begin{eqnarray*}
-\intt_0^{2\pi}\log|\g(re^{i\theta})|\frac{d\theta}{2\pi} &\ge& m\log\l(\frac{\sqrt{m}}{r}\r) - \l(\frac{1}{2}+\eps\r)m \\
&=& \frac{1}{2} m\log(m) - m\log(r) -\l(\frac{1}{2}+\eps\r)m \\
&=& \frac{1}{2} m\log(m) - O(m)
\end{eqnarray*}
 Thus
\begin{eqnarray*}
\P[n(r)\ge m] &\le& \P\l[\log M(R_m)\ge \l(\frac{1}{2}+\eps\r)m\r] + \P\l[- \intt_0^{2\pi}\log|\g(re^{i\theta})|\frac{d\theta}{2\pi}\ge \frac{1}{2} m\log(m)-O(m) \r]\\
 &\le& e^{-e^{\eps m}} + \P\l[-
 \intt_0^{2\pi}\log|\g(re^{i\theta})|\frac{d\theta}{2\pi}\ge
 \frac{1}{2}m\log(m)(1+o(1)) \r] \hsp{1.2cm} \mb{by }\ref{eq:jensens}. 
\end{eqnarray*}
From Lemma~\ref{lem:ubdonintegral}, we deduce that for any $\del>0$, there is a constant $C_2$ such that  
\begin{eqnarray*} 
\P\l[- \intt_0^{2\pi}\log|\g(re^{i\theta})|\frac{d\theta}{2\pi} \ge \frac{1}{2} m\log(m) (1+o(1))\r] &\le& C_2 e^{-(2-\del)(\frac{m}{2} \log(m))^2/\log(\frac{m}{2} \log(m))}\\
&\le& C_2 e^{-(\frac{1}{2}-\frac{\del}{4})m^2\log(m) (1+o(1))}. 
\end{eqnarray*} From this, the upper bound follows. 
\end{proof}

\begin{lemma}\label{lem:ubdonintegral} For any given $\del>0$, $\exists C_2$ such that $\P\l[- \intt_0^{2\pi}\log|\g(re^{i\theta})|\frac{d\theta}{2\pi}\ge m \r]\le C_2e^{-\frac{(2-\del)m^2}{\log(m)}}$ $\forall m$.

\end{lemma}
\begin{proof} Let $P$ be the Poisson kernel on $D(0;r)$. Fix $\eps>0$
  and  let
 $A_{\eps}=\sup\{P(re^{i\theta},w):|w|=\eps,\theta\in[0,2\pi) \}$ and
 $B_{\eps}=\inf\{P(re^{i\theta},w):|w|=\eps,\theta\in[0,2\pi) \}$.
Since $\log |\g|$ is a sub-harmonic function, for any $w$ with $|w|=\eps$, we get
\begin{eqnarray*}
\log |\g(w)| &\le& \intt_0^{2\pi}\log|\g(re^{i\theta})|P(re^{i\theta},w)\frac{d\theta}{2\pi} \\
&\le& A_{\eps} \intt_0^{2\pi}\log_+|\g(re^{i\theta})|\frac{d\theta}{2\pi} - B_{\eps}\intt_0^{2\pi}\log_- |\g(re^{i\theta})|\frac{d\theta}{2\pi} \\
&=&  A_{\eps} \intt_0^{2\pi}\log_+|\g(re^{i\theta})|\frac{d\theta}{2\pi} + B_{\eps} \l( \intt_0^{2\pi}\log|\g(re^{i\theta})|\frac{d\theta}{2\pi}- \intt_0^{2\pi}\log_+|\g(re^{i\theta})|\frac{d\theta}{2\pi}\r) \\
&\le& B_{\eps} \intt_0^{2\pi}\log|\g(re^{i\theta})|\frac{d\theta}{2\pi} + A_{\eps} \log_+ M(r).
\end{eqnarray*}
This implies $\log M(\eps) \le B_{\eps}
\intt_0^{2\pi}\log|\g(re^{i\theta})|\frac{d\theta}{2\pi} + A_{\eps}
\log_+ M(r)$.

Therefore if $\intt_0^{2\pi}\log|\g(re^{i\theta})|\frac{d\theta}{2\pi} \le -m$, then one of the following must happen. Either $\{\log M(\eps)\le -B_{\eps} m+\sqrt{m}$\} or  $\{A_{\eps}\log_+M(r)>\sqrt{m}\}$. 

Using (\ref{eq:ubdformax}), since $M(r)<M\l(Cm^{\frac{1}{4}}\r)$ for
any $C$, we see that $\P\l[\log_+M(r) >\frac{\sqrt{m}}{A_{\eps}}\r]\le
e^{-e^{Cm}}$ for some  constant $C$ depending on $\eps$. Hence
\begin{eqnarray*}
\P\l[\intt_0^{2\pi}\log|\g(re^{i\theta})|\frac{d\theta}{2\pi} \le -m \r] &\le& e^{-e^{Cm}} + \P\l[ \log M(\eps) \le -B_{\eps} m+\sqrt{m} \r] \\
 &\le& e^{-e^{Cm}} + e^{-2B_{\eps}^2\frac{m^2}{\log(m)}(1+o(1))}
\end{eqnarray*}
where in the last line we have used Lemma~\ref{lem:maxmodulus}.

As $\eps\tends 0$, $B_{\eps}\tends 1$ and hence the proof is complete.

\end{proof}

Now we prove the upper bound on the maximum modulus in a disk of radius $r$ that was used in the last part of the proof of Lemma~\ref{lem:ubdonintegral}. For possible future use we prove a lower bound too.

\begin{lemma}\label{lem:maxmodulus} Fix $r>0$. There are constants $\alp,C_1,C_2$ such that 
\[ C_1 e^{-\frac{\alp m^2}{\log(m)}} \le \P[\log M(r)\le -m] \le C_2 e^{-\frac{2m^2}{\log(m)}(1+o(1))}.\]
\end{lemma}
\begin{proof}

{\bf Lower bound } By Cauchy-Schwarz, 
$M(r)\le \l(\summ_{n=0}^{k-1} |a_n|^2 \r)^{1/2}e^{r^2/2} +
\summ_{n=k}^{\infty} \frac{|a_n|r^n}{\sqrt{n!}}$. We shall choose $k$
later. We will bound from below the probability that each of these
summands is less than $\frac{e^{-m}}{2}$.  

Let $\phi_k$ denote the density of $\Gam(k,1)$.
\begin{eqnarray*}
\P\l[\l(\summ_{n=0}^{k-1} |a_n|^2 \r)^{1/2}e^{r^2/2}\le
\frac{e^{-m}}{2}\r] &=&\P\l[\summ_{n=0}^{k-1} |a_n|^2 \le
\frac{e^{-2m}e^{r^2}}{4} \r]\\ 
&\ge& \phi_k\l(\frac{e^{-2m}e^{r^2}}{8}\r)\frac{e^{-2m}e^{r^2}}{8} \\
&=& e^{-2mk-k\log(k)+O(k)} 
\end{eqnarray*}

Also if $|a_n|\le n^2$ $\forall n\ge k$, then the second summand
\[ \summ_{n=k}^{\infty} |a_n|\frac{r^n}{\sqrt{n!}} \le C\frac{r^kk^2}{\sqrt{k!}} \le C e^{-k\log(k)/3} \]
Also the event $\{|a_n|\le n^2$ $\forall n\ge k \}$ has probability at least $1-\summ_{n=k+1}^{\infty} e^{-n^4} \ge 1-Ce^{-k^4}$.

Thus if we set $k=\frac{\gam m}{\log(m)}$ for a sufficiently large $\gam$, then both the terms are less than $e^{-\frac{m}{2}}$ with probability at least $e^{-2\gam m^2/\log(m)}$. 

{\bf Upper bound } By Cauchy's theorem, \[ a_n=\frac{\sqrt{n!}}{2\pi i}\intt_{C_r} \frac{\g(\zet)}{\zet^{n+1}} d\zet, \] where $C_r$ is the curve $C_r(t)=re^{it}$, $0\le t\le 2\pi$. Therefore, 
\[ |a_n| \le \frac{M(r)\sqrt{n!}}{r^n}. \]
Thus we get
\[ \P[M(r)\le e^{-m}] \le \prodd_{n=0}^{\infty} \P\l[|a_n|\le \frac{e^{-m}\sqrt{n!}}{r^n}\r]. \]
$|a_n|^2$ are i.i.d.  exponential random variables with mean $1$. Therefore, $\P\l[|a_n|\le \frac{e^{-m}\sqrt{n!}}{r^n}\r]\le  \frac{e^{-2m} n!}{r^{2n}}$. Using this bound for $n\le k :=\frac{\bet m}{\log(m)}$, we get
\begin{eqnarray*}
\P[M(r)\le e^{-m}] &\le& \prodd_{n=0}^k  \frac{e^{-2m} n!}{r^{2n}} \\
&\le& C e^{-2mk + \frac{k^2}{2}\log(k)+O(k^2)} \\
&\le& Ce^{(-2\bet+\frac{\bet^2}{2})\frac{m^2}{\log(m)}+O(\frac{m^2}{(\log(m))^2})}.
\end{eqnarray*}
$-2\bet+\frac{\bet^2}{2}$ is minimized when $\bet=2$ and we get,
\begin{equation}
 \P[M(r)\le e^{-m}] \le e^{-2\frac{m^2}{\log(m)}(1+o(1))}.  
\end{equation}
\end{proof}

\section{Overcrowding - The Hyperbolic case}\label{sec:hyperbolic}

\subsection{Case $\rho=1$}
We give a quick proof of Theorem~\ref{thm:hyperbolic} in the special
case $\rho=1$, as it is much easier and moreover we get matching upper
and lower bounds. The proof is similar to the case of the Ginibre
ensemble dealt with in Theorem~\ref{thm:ginibre} and is based on the
fact that the set of absolute values of the zeroes of $\f_1$ is
distributed the same as a certain set of independent random
variables. The reason for this similarity between the two cases owes
to the fact that both of them are determinantal. The zero set of
$\f_1$ is a determinantal process with the 
Bergman kernel for the unit disk, namely 
\[ K_B(z,w)=\frac{1}{\pi}\frac{1}{(1-z{\bar w})^2}, \]
  as discovered by Peres and Vir{\'a}g~\cite{pervir}. 

{\bf Proof of Theorem~\ref{thm:hyperbolic} for $\rho=1$} 

 By Peres and Vir{\'a}g~\cite{pervir}, Theorem 2 (ii), the set of
 absolute values of the zeroes of $\f_1$ has the same distribution as
 the set $\{U_n^{1/2n}\}$ where $U_n$ are i.i.d.  Uniform$[0,1]$ random
 variables. Therefore, 
\begin{eqnarray*} 
\P\l[n(r)\ge m\r] &\ge& \prodd_{n=1}^{m} \P\l[U_n^{1/2n}<r \r] \\
 &=& \prodd_{n=1}^{m} r^{2n} \\
 &=& r^{m(m+1)}.
\end{eqnarray*}
To prove the inequality in the other direction, note that
\begin{eqnarray*}
\P\l[n(r)\ge m\r] &\le& \P\l[\summ_{n=1}^{m^2}\ind(U_n^{1/2n}<r)\ge m \r]+\summ_{n=m^2+1}^{\infty} \P\l[U_n^{1/2n}<r\r] \\
&\le& {m^2 \choose m}\prodd_{n=1}^m\P\l[U_n^{1/2n}<r \r] + \summ_{n>m^2}r^{2n}\\
&=& {m^2 \choose m}r^{m(m+1)}+ \frac{r^{2m^2+2}}{1-r^2}\\
&=& r^{m(m+1)}\l(1+ O\l(e^{m\log(m)}\r)\r).
\end{eqnarray*}
This completes the proof of the theorem for $\rho=1$. 

\subsection{All values of $\rho$}


%



{\bf Remark: } Overall, the idea of proof is the same as that of Theorem~\ref{thm:planar}. However we do not get matching upper and lower bounds in the present case, the reason being that in the hyperbolic analogue of Lemma~\ref{lem:maxmodulus}, the leading term in the exponent of the upper bound does depend on $r$, unlike in the planar case. (An examination of the proof of Theorem~\ref{thm:planar} reveals that we get a matching upper bound only because replacing $r$ by $\eps$  does not affect the leading term in the exponent in the upper bound in Lemma~\ref{lem:maxmodulus}). However we still expect that the lower bound in Theorem~\ref{thm:hyperbolic} is tight. (See remark after the proof).

\begin{proof}[Proof of Theorem~\ref{thm:hyperbolic}] 

{\bf Lower Bound }  As before we find a lower bound for the probability that the $m^{\mb{th}}$ term dominates the rest. Note that if $|z|=r$, 
\begin{equation}
\Mid\f_{\rho}(z)-{-\rho \choose m}^{1/2} a_mz^m \Mid \le \summ_{n=0}^{m-1} |a_n| {-\rho \choose n}^{1/2}r^n + \summ_{n=m+1}^{\infty} |a_n| {-\rho \choose n}^{1/2}r^n 
\end{equation}\label{eq:threeconditions}
Now suppose the following happen-
\begin{enumerate}
\item $|a_n| \le \sqrt{n}$ $\forall n\ge m+1$.
\item $|a_m| \ge (\alp+1)\sqrt{m}$ where $\alp$ will be chosen shortly.
\item $|a_n|{-\rho\choose n}^{1/2}r^n < \frac{1}{\sqrt{m}}{-\rho\choose m}^{1/2}r^m$ for every $0\le n\le m-1$.
\end{enumerate}
Then the right hand side of (\ref{eq:threeconditions}) is bounded by
\begin{eqnarray*}
\mb{RHS of } (\ref{eq:threeconditions}) &\le& \summ_{n=0}^{m-1}|a_n|{-\rho\choose n}^{1/2}r^n  + \summ_{n=m+1}^{\infty}|a_n|{-\rho\choose n}^{1/2}r^n \\
&\le& \summ_{n=0}^{m-1}\frac{1}{\sqrt{m}} {-\rho\choose m}^{1/2}r^m + \summ_{n=m+1}^{\infty}\sqrt{n}{-\rho\choose n}^{1/2}r^n \\
&\le& \sqrt{m} {-\rho\choose m}^{1/2}r^m + C \sqrt{m}{-\rho\choose m}^{1/2}r^m \hsp{1.5cm} \mbox{ for some }C \\
&=&  (C+1)\sqrt{m} {-\rho\choose m}^{1/2}r^m\\
&\le& |a_m| {-\rho\choose m}^{1/2}r^m
\end{eqnarray*}
if $\alp=C$. Thus if the above three events occur with $\alp=C$, then the $m^{\mb{th}}$ term dominates the sum of all the other terms on $\d D(0;r)$. Also these events have probabilities as follows. 
\begin{enumerate}
\item $\P[|a_n| \le \sqrt{n}$ $\forall n\ge m+1]\ge 1- \summ_{n=m+1}^{\infty} e^{-n}\ge 1-C'e^{-m}$. 
\item $\P[|a_m|\ge (\alp+1)\sqrt{m}] = e^{-(\alp+1)^2m}$.
\item The third event has probability as follows. Recall again that
  $\P\l[\xi<x\r]\ge \frac{x}{2}$ if $x<1$ and $\xi$ is Exponential with
  mean $1$. We apply this below with $x=\l(\frac{{-\rho\choose
  m}^{1/2}r^{m-n}}{\sqrt{m} {-\rho\choose n}^{1/2}}\r)^2$. This is
  clearly less than $1$. Thus 
\begin{eqnarray*}
\P\l[|a_n|\le \frac{ {-\rho\choose m}^{1/2}r^{m-n}}{\sqrt{m} {-\rho\choose n}^{1/2} } \hsp{2mm}\forall n\le m-1 \r] &=& \prodd_{n=0}^{m-1} \P\l[|a_n|\le \frac{ {-\rho\choose m}^{1/2}r^{m-n}}{\sqrt{m} {-\rho\choose n}^{1/2} } \r] \\
&\ge& \prodd_{n=0}^{m-1}\frac{{-\rho\choose m}r^{2m-2n}}{2 m {-\rho\choose n} } \\
&=& r^{m(m+1)} m^{-m} \prodd_{n=0}^{m-1} \frac{(m+1)\ldots (m+\rho-1)}{ (n+1)\ldots (n+\rho-1)} \\
&\ge& r^{m(m+1)} m^{-m} \prodd_{n=0}^{m-1}\frac{m^{\rho}}{(n+\rho)^{\rho}} \\
&\ge& r^{m(m+1)+O(m\log(m))}.
\end{eqnarray*}
\end{enumerate}
Since these three events are independent, we get the lower bound in the theorem. 

{\bf Upper Bound } The proof will proceed along the same lines as in
Theorem~\ref{thm:planar}. We need the following analogue of
Lemma~\ref{lem:maxmodulus}.  

\begin{lemma}\label{lem:hypmaxmod} Fix $r<1$.  Let $M(r)=\sup_{z\in D(0;r)} |\f_{\rho}(z)|$. Then
\[ \P[M(r)\le e^{-m}] \le e^{-\frac{m^2}{|\log(r)|}(1+o(1))}. \]
\end{lemma}
\begin{proof}
By Cauchy's theorem, for every $n\ge 0$,
\[ a_n{-\rho \choose n}^{1/2} = \frac{1}{2\pi i} \intt_{rT} \frac{\f(\zet)}{\zet^{n+1}} d\zet. \]
From this we get
\[ |a_n|^2 \le \frac{M(r)^2}{{-\rho \choose n} r^{2n}}. \]
Since ${-\rho \choose n}\ge \frac{n^{\rho-1}}{\Gam(\rho+1)}$, we obtain
\begin{eqnarray*}
\P[M(r)\le m] &\le& \prodd_n \P[|a_n|^2\le \frac{\Gam(\rho+1)e^{-2m}}{n^{\rho-1}r^{2n}}] \\
&\le& \prodd_{n=0}^\frac{m}{\log(1/r)}\frac{\Gam(\rho+1) e^{-2m}}{r^{2n}n^{\rho-1}} \\
&\le& e^{-\frac{2m^2}{\log(1/r)}+\l(\frac{m}{\log(1/r)}\r)^2\log(r)+O(m\log(m))}\\
&=& e^{-\frac{m^2}{\log(1/r)}+O(m\log(m))}. 
\end{eqnarray*}
\end{proof}

Coming back to the proof of the upper bound in the theorem, fix $R$ such that $r<R<1$. Then by Jensen's formula,
\begin{equation}\label{eq:jensenss}
n(r)\log\l(\frac{R}{r}\r) \le \intt_r^R \frac{n(u)}{u} du =
\intt_{R\T}\log|f(Re^{i\theta})| \frac{d\theta}{2\pi} - \intt_{r\T}
\log|f(re^{i\theta})| \frac{d\theta}{2\pi}.\end{equation}

Now consider the first summand in the right hand side of (\ref{eq:jensens}).
\begin{eqnarray*}
\P\l[\intt_{R\T}\log|\f(Re^{i\theta})|\frac{d\theta}{2\pi}
>\sqrt{m}\r] &\le& \P\l[\log M(R)\ge \sqrt{m}\r].
\end{eqnarray*}
Now suppose that $|a_n|<\lam^n$ $\forall n\ge m+1$ where
$1<\lam<1/R$. This has probability at least
$C_1e^{-\lam^{2m}/2}$. Then, 
\begin{eqnarray*}
M(R) &\le& \summ_{n-0}^{\infty} |a_n|{-\rho \choose n}^{1/2} R^n \\
&\le& \l(\summ_{n=0}^m |a_n|^2 \r)^{1/2}C_R + C_{R'}
\end{eqnarray*}
for some constants $C_R$ and $C_{R'}$.

Thus if $M(R)>e^{\sqrt{m}}$  then either $\summ_{n=0}^m |a_n|^2 >Ce^{2\sqrt{m}}$ or else $|a_n|>\lam^n$ for some  $n\ge m+1$. Thus
\begin{equation*}
\P\l[M(R)>\sqrt{m}\r] \le e^{-e^{c\sqrt{m}}}.
\end{equation*}
This proves that
\begin{equation*}
\P\l[\intt_{R\T}\log|\f(Re^{i\theta})|\frac{d\theta}{2\pi} >\sqrt{m}\r]\le e^{-e^{c\sqrt{m}}}.
\end{equation*}
Fix $\delta>0$ and $R$ close enough to $1$ such that
$\log(R)>-\delta$. Then with probability $\ge 1-e^{-e^{c\sqrt{m}}}$,
we obtain from (\ref{eq:jensenss}), 
\[ -\intt_{r\T} \log|\f(re^{i\theta})|\frac{d\theta}{2\pi} \ge m\l(\log\l(\frac{1}{r}\r) -\delta \r)-\sqrt{m}. \]

Now the calculations in the proof of Lemma~\ref{lem:ubdonintegral} show that
\[ \log M(\eps) \le B_{\eps} \intt_0^{2\pi} \log|\f(re^{i\theta})| P(re^{i\theta},w) \frac{d\theta}{2\pi}  + A_{\eps}\log_+ M(r). \]
Here $0<\eps<r$ is arbitrary and $A_{\eps},B_{\eps}$ are as defined in Lemma~\ref{lem:ubdonintegral}.  By the same computations as in that Lemma, we obtain,
we obtain the inequality
\begin{equation*} \P\l[\intt_0^{2\pi} \log |\f(re^{i\theta})|\frac{d\theta}{2\pi}\le
-m(|\log r|-\del)+\sqrt{m} \r] \le
 e^{-B_{\eps}^2\frac{m^2\log^2(r)(1-\delta)}{|\log (\eps)|}}+e^{-e^{cm}}.
\end{equation*}
Therefore, by (\ref{eq:jensenss})
\begin{equation*}
  \P\l[n(r)\ge m \r] \le e^{-\kappa m^2\log^2(r)(1+o(1))},
\end{equation*}
where $\kappa =\sup \l\{\frac{B_{\eps}^2}{|\log (\eps)|}: 0<\eps <r \r\}$.
 However it is clear that this cannot be made to match the lower bound
 by any choice of $\eps$.  
\end{proof}
 
{\bf Remark :} If we could prove   
\begin{equation*}\P\l[\intt_0^{2\pi} \log
  |\f(re^{i\theta})|\frac{d\theta}{2\pi}\le -x \r] \le
  e^{-\frac{x^2}{|\log(r)|}}, 
\end{equation*} 
that would have given us a matching upper bound. Now, one way for the
event $\intt_0^{2\pi} \log |\f(re^{i\theta})|\frac{d\theta}{2\pi}\le
-x$ to occur is to have $\log M(r)<-x$ which,  by
Lemma~\ref{lem:hypmaxmod} has probability at most
$e^{-x^2/\log(\frac{1}{r})}$. One way to proceed could be to show that
if the integral is smaller than $-x$, so is $\log M(s)$ for $s$
arbitrarily close to $r$ (with high probability). Alternately, if we
could bound the coefficients directly by the bound on the integral (as
in Lemma~\ref{lem:hypmaxmod}), that would also give us the desired
bound. For these reasons, and keeping in mind the case $\rho=1$, where we do have a
matching upper bound, we believe that the lower bound in Theorem~\ref{thm:hyperbolic} is tight.

\section{Moderate and Very Large deviations for the planar GAF}\label{sec:moderate}
In this section we prove Theorem~\ref{thm:verylarge} and Theorem~\ref{thm:moderate}.
\begin{remark}In the case $\alp\ge 2$, one side of the estimate as asked for in
  the conjecture (with $\log \log$ of the probability) follows
  trivially from the results in  
Sodin and Tsirelson~\cite{sodtsi3}. They prove that for any $\del>0$,
there exists a constant $c(\del)$ such that 
\[ \P\l[ |n(r)-r^2| >\del r^2 \r] \le e^{-c(\del)r^4}. \]
When $\alp\ge 2$, clearly $n((1-\del)r^{\sqrt{\alp}})\ge n(r)$, whence
from the above result it follows that  
\begin{eqnarray*}
\P\l[n(r)\ge r^2+r^{\alp} \r] &\le& \P\l[n((1-\del)r^{\sqrt{\alp}})\ge r^{\alp} \r]\\
&\le& e^{-c(\del)r^{2\alp}}.
\end{eqnarray*}
This gives 
\begin{equation}\label{eq:mildbd}
\limsup_{r\tends \infty} \frac{\log \log
         \l(\frac{1}{\P[|n(r)-r^2|>r^{\alp}]}\r)}{\log r} \le  2\alp.
\end{equation}
\end{remark}

The obviously loose inequality $n((1-\del)r^{\sqrt{\alp}})\ge n(r)$
that we used, suggests that (\ref{eq:mildbd}) can be improved when $\alp>2$
to Theorem~\ref{thm:verylarge}. 

\begin{proof}[Proof of Theorem~\ref{thm:verylarge}]

{\bf Lower Bound }Let $m=r^2+\gam r^{\alp}$. Suppose the $m^{\mb{th}}$ term dominates the sum of all the  other terms on $\d D(0;r)$, i.e., suppose
\begin{equation}\label{eq:dominate2} 
\Mid\frac{a_m z^m }{\sqrt{m!}}\Mid \ge \Mid \summ_{n\neq m} \frac{a_n z^n}{\sqrt{n!}}\Mid \hspace{2cm} \mb{ whenever }|z|=r.
\end{equation}

Now we want to find a lower bound for the probability of the event in (\ref{eq:dominate2}). Note that the left side of (\ref{eq:dominate2}) is identically equal to $\frac{|a_m|r^m}{\sqrt{m!}}$. 

Now suppose the following happen-
\begin{enumerate}
\item $|a_n| \le n$ $\forall n\ge m+1$.
\item $|a_m| \ge m$.
\item $|a_n|\frac{r^{n}}{\sqrt{n!}} < \frac{\gam r^{\alp}}{m}\frac{r^m}{\sqrt{m!}}$ for every $0\le n\le m-1$.
\end{enumerate}
Then the right hand side of (\ref{eq:dominate2}) is bounded by
\begin{eqnarray*}
\mb{RHS of } (\ref{eq:dominate2}) &\le& \summ_{n=0}^{m-1}|a_n|\frac{r^{n}}{\sqrt{n!}}  + \summ_{n=m+1}^{\infty} \frac{|a_n|r^n}{\sqrt{n!}} \\
&\le& \summ_{n=0}^{m-1}\frac{\gam r^{\alp}}{m}\frac{r^{m}}{ \sqrt{m!}} + \summ_{n=m+1}^{\infty} \frac{nr^n}{\sqrt{n!}} \\
&\le& \frac{mr^{m}}{\sqrt{m!}}\l( \frac{\gam r^{\alp}}{ m} + o(1) \r) \\
&\le& |a_m|\frac{r^m}{m!}
\end{eqnarray*}

 Thus if the above three events occur, then the $m^{\mb{th}}$ term dominates the sum of all the other terms on $\d D(0;r)$. Also these events have probabilities as follows.
\begin{enumerate}
\item $\P[|a_n| \le n$ $\forall n\ge m+1]\ge 1- \summ_{n=m+1}^{\infty} e^{-n^2}\ge 1-C'e^{-m^2}=1-o(1)$. 
\item $\P[|a_m|\ge m] = e^{-m^2}=e^{-\gam^2 r^{2\alp}(1+o(1))}$.
\item The third event has probability as follows. Recall again that
  $\P\l[\xi<x\r]\ge \frac{x}{2}$ if $x<1$ and $\xi$ is Exponential with
  mean $1$. We apply this below with $x=\frac{\gam^2 r^{2\alp}}{m^2}\frac{r^{2m-2n}n!}{
    m!}$. This is clearly less than $1$ if $n\ge
  r^2$. Therefore if $m$ is sufficiently large it is easy to see that
  for all $0\le n\le m-1$, the same is valid. Thus
\begin{eqnarray*}
\P\l[|a_n|\le \frac{\gam r^{\alp}}{m}\frac{r^{m-n}\sqrt{n!}}{\sqrt{m!}} \hsp{2mm}\forall n\le m-1 \r] &=& \prodd_{n=0}^{m-1} \P\l[|a_n|\le \frac{\gam r^{\alp}}{m}\frac{r^{m-n}\sqrt{n!}}{\sqrt{m!}} \r] \\
&\ge& \prodd_{n=0}^{m-1} \frac{\gam^2 r^{2\alp}}{m^2}\frac{r^{2m-2n}n!}{2 m!} \\
&=& r^{2\alp (m+1)+m(m+1)} 2^{-m}m^{-2m}e^{-\summ_{k=1}^m k\log k}  \\
&=& e^{m^2\log(r)-\frac{1}{2}m^2\log(m)+O(m^2)}\\
&=& e^{-(\frac{\alp}{ 2}-1)\gam^2 r^{2\alp} \log(r) +O(r^{2\alp})}
\end{eqnarray*}

\end{enumerate}

Since these three events are independent, we get 
\begin{equation}
\P\l[n(r)\ge r^2+r^{\alp} \r] \ge e^{-\l(\frac{\alp}{ 2}-1\r)\gam^2 r^{2\alp}\log r + O(r^{2\alp})}.
\end{equation}
{\bf Upper Bound} We omit the proof of the upper bound, as it follows 
the same lines as that of Theorem~\ref{thm:planar} and we have
already seen such arguments again in the proof of
Theorem~\ref{thm:hyperbolic} (In those two cases as well as the
present case, we are looking at very large deviations, and that is the
reason why the same tricks work). 

Moreover note that the lower bound along with (\ref{eq:mildbd}) proves the
statement in the conjecture. 
\end{proof}

{\bf Case $1<\alp<2$: }We prove Theorem~\ref{thm:moderate}. Along with
Theorem~\ref{thm:verylarge} this shows that the asymptotics of $\P\l[
n(r)\ge r^2+\gam r^{\alp}\r]$ does undergo a qualitative change at
$\alp=2$.

\begin{proof}[Proof of Theorem~\ref{thm:moderate}] Write $m=r^2+\gam
  r^{\alp}$. As usual, we bound $\P\l[n(r)\ge m\r]$ from below by the
  probability of the event that the $m^{\mb{th}}$ term dominates the
  rest of the series.

Firstly, we need a couple of estimates. Consider $\frac{r^{2n}}{n!}$
as a function of $n$. This increases 
monotonically  up to $n=r^2$ and then decreases
monotonically. $m=r^2+\gam r^{\alp}$ is on the latter part. Write
$M=r^2-\gam r^{\alp}$.

Firstly,  observe that $(r^2-k)(r^2+k)<(r^2)^2$, for  $1\le k\le
\gam r^{\alp}$, whence $r^{2m-2M}>\prodd_{j=M+1}^{m-1} j$. This
implies that 
\begin{equation}\label{eq:relationMm}
\frac{r^M}{\sqrt{M!}}<\frac{r^m}{\sqrt{m!}}
\end{equation}

Secondly, note that for any $n=M-p$, 
\begin{eqnarray*}
  \frac{r^{2n}/n!}{r^{2M}/M!} &=&\prodd_{j=0}^{p-1} \frac{M-j}{r^2} \\
             &=& \prodd_{j=0}^{p-1} (1-\gam r^{\alp-2}-jr^{-2}) \\
             &\le& e^{-\summ_{j=0}^{p-1} (\gam r^{\alp-2}+jr^{-2})} \\
             &=& e^{-\gam pr^{\alp-2} -\frac{p(p+1)}{2}r^{-2}}.
\end{eqnarray*}
Now we set $p=Cr^{2-\alp}$ with $C$ so large that 
$ e^{-\gam C}\le \frac{1}{4}$.

Then also note that if $n<M-kp$, it follows that
\begin{equation}
  \label{eq:decayofcoeffs}
  \frac{r^{2n}/n!}{r^{2m}/m!}\le \frac{1}{4^k},
\end{equation}
where we used (\ref{eq:relationMm}) to replace $M$ by $m$.

Thirdly, if $n=m+p$ with $p\le r^2-\gam r^{\alp}$, then,
\begin{eqnarray*}
  \frac{r^{2n}/n!}{r^{2m}/m!} &=&\prodd_{j=1}^{p} \frac{r^2}{m+j} \\
             &=& \prodd_{j=1}^{p} (1+\gam r^{\alp-2}+jr^{-2})^{-1} \\
             &\le& e^{-\frac{1}{2}\summ_{j=1}^{p} (\gam
             r^{\alp-2}+jr^{-2})} \\ 
             &=& e^{-\frac{1}{2}(\gam p r^{\alp-2} +\frac{p(p+1)}{2}r^{-2})}.
\end{eqnarray*}
If $p=2Cr^{2-\alp}$, where $C$ was as chosen before, then for
$n>m+kp$, we get   
\begin{equation}
  \label{eq:decayofcoeffs2}
  \frac{r^{2n}/n!}{r^{2m}/m!}\le \frac{1}{4^k}.
\end{equation}
From now on $p=2Cr^{2-\alp}$ is fixed so that (\ref{eq:decayofcoeffs})
and (\ref{eq:decayofcoeffs2}) are satisfied.

Next we divide the coefficients other than $m$ into groups: 
\begin{itemize}
\item $A_k=\{n:n\in (M-kp,M-(k-1)p]$ for $1\le k\le \lceil \frac{M}{p}
  \rceil \}$.
\item $D_k=\{n:n\in [m+(k-1)p,m+kp)$ for $1\le k\le \lceil \frac{M}{p}
  \rceil$\}.
\item $B=\{n:n\in [M+1,m-1]\}$.
\item $C=\{n:n\in [2r^2,\infty)\}$.
\end{itemize}
\begin{remark}
  As defined, there is an overlap between $D_{\lceil \frac{M}{p}
  \rceil}$ and $C$. This is inconsequential, but for definiteness, let
  us truncate the former interval at $r^2$ (just as $A_{\lceil \frac{M}{p}
  \rceil}$ is understood to be truncated at $0$).
\end{remark}

Now consider the following events.

\begin{enumerate}
\item $|a_n| \le \frac{2^k}{M}$ for $n\in A_k$ for $k\le \lceil
  \frac{M}{p}\rceil$ \}.
\item $|a_n|\le \frac{2^k}{M}$ for $n\in D_k$ for $k\le
  \lceil \frac{M}{p}\rceil \}$.
\item $\summ_{n\in B} |a_n|\frac{r^n}{\sqrt{n!}} \le 4 \frac{r^m}{\sqrt{m!}}$.
\item $|a_n|<n-2r^2$ for $n\in C$.
\item $|a_m| \ge 15$.
\end{enumerate}

Suppose all these events occur. Then

\begin{enumerate}
\item The event $|a_n| \le \frac{2^k}{M}$ for $n\in A_k$, $k\le \lceil
  \frac{M}{p}\rceil$ gives
\begin{eqnarray}\label{eq:A}
  \sup\{\Mid \summ_{n=0}^M \frac{a_nz^n}{\sqrt{n!}} \Mid:|z|=r \} &\le&
  \summ_{k=1}^{\lceil M/p \rceil} \summ_{n=M-kp+1}^{M-(k-1)p}
  |a_n|\frac{r^n}{\sqrt{n!}} \\
     &\le&  \summ_{k=1}^{\lceil M/p
  \rceil}\frac{1}{2^{k}}\frac{r^{m}}{\sqrt{m!}}\frac{2^kp}{M}
  \hsp{1cm} \mb{by }(\ref{eq:decayofcoeffs})\\
 &\le& \frac{r^{m}}{\sqrt{m!}} \summ_{k=1}^{\lceil M/p \rceil}\frac{p}{M} \\
   &\le& \frac{r^{m}}{\sqrt{m!}} \l(1+\frac{p}{M}\r). 
\end{eqnarray}

\item The event $|a_n|\le \frac{2^k}{M}$ for $n\in D_k$, $k\le
  \lceil \frac{M}{p}\rceil$ gives
\begin{eqnarray}\label{eq:D}
\sup\{\Mid \summ_{n=m+1}^{2r^2} \frac{a_nz^n}{\sqrt{n!}}\Mid:|z|=r \} &=& 
\summ_{k=1}^{\lceil M/p \rceil} \summ_{n=m+(k-1)p+1}^{M+kp}
  |a_n|\frac{r^n}{\sqrt{n!}} \\
 &\le&\summ_{k=1}^{\lceil M/p
  \rceil}\frac{1}{2^{k}}\frac{r^{m}}{\sqrt{m!}}\frac{2^kp}{M}
  \hsp{1cm} \mb{by }(\ref{eq:decayofcoeffs2}) \\
 &\le& \frac{r^m}{\sqrt{m!}} \l(1+\frac{p}{M}\r). 
\end{eqnarray}

\item The third event gives 
\begin{equation}\label{eq:B}
\summ_{n\in B} |a_n|\frac{r^n}{\sqrt{n!}} \le 4 \frac{r^m}{\sqrt{m!}},
\end{equation}
by assumption.

\item The event $|a_n|<n-2r^2$ for $n\in C$: Since $n>2r^2$,
\begin{eqnarray*}
\frac{r^n}{\sqrt{n!}} &=& \frac{r^m}{\sqrt{m!}}\prodd_{k=m+1}^{n} \frac{r}{\sqrt{k}} \\
&\le&  \frac{r^m}{\sqrt{m!}}\prodd_{k=2r^2+1}^{n} \frac{r}{\sqrt{k}} \\
&\le& \frac{r^m}{\sqrt{m!}}\l(\frac{1}{\sqrt{2}}\r)^{n-2r^2}.
\end{eqnarray*}
Therefore we get (using $|a_n|<n-2r^2$ $\forall n>2r^2$)
\begin{eqnarray}
\summ_{n\in C} |a_n|\frac{r^n}{\sqrt{n!}} &\le& \frac{r^m}{\sqrt{m!}} \summ_{n>2r^2} (n-2r^2) \l(\frac{1}{\sqrt{2}}\r)^{n-2r^2} \\
 &=& \frac{r^m}{\sqrt{m!}} \frac{\sqrt{2}}{(\sqrt{2}-1)^2}.
\end{eqnarray}
\end{enumerate}

Putting together the contributions from these four groups of terms,
and using $|a_m|>15$, we get (for large values of $r$)
\begin{eqnarray*}
\summ_{n\not=m}|a_n|\frac{r^n}{\sqrt{n!}} \le |a_m|\frac{r^m}{\sqrt{m!}}.
\end{eqnarray*}

Now we compute the probabilities of the events enumerated above.

\begin{enumerate}
\item The event $|a_n| \le \frac{2^k}{M}$ for $n\in A_k$ for $k\le \lceil
  \frac{M}{p}\rceil$.
Now for a fixed $k\le 3\log_2(r)$, we deduce
\begin{eqnarray*}
  \P\l[|a_n|\le \frac{2^k}{M} \mb{ for } n\in A_k \r] &\ge&
  \P\l[|a_n|\le \frac{1}{M}  \mb{ for } n\in A_k \r] \\
   &\ge& \l(  \frac{1}{2M^2} \r)^{p}.
\end{eqnarray*}
Therefore
\begin{eqnarray*}
  \P\l[|a_n|\le \frac{2^k}{M} \mb{ for } n\in A_k \mb{
  for every }k\le 3\log_2(r)\r] &\ge& \l(  \frac{1}{2M^2}
  \r)^{3p\log_2(r)}\\
&\ge& e^{-cr^{2-\alp}(\log(r))^2}.
\end{eqnarray*}
for some $c$.

Next we deal with $k> 3\log_2(r)$.

\begin{eqnarray*}
  \P\l[|a_n|\le \frac{2^k}{M} \mb{ for } n\in A_k \mb{
  for every }k>3\log_2(r)\r] &\ge& 1- \summ_{k>3\log_2(r)}p \P\l[|a|>
  \frac{2^k}{M}\r]  \\
 &=& 1- \summ_{k>3\log_2(r)}p e^{-2^{2k}M^{-2}}.
\end{eqnarray*}
Now the summation in the the last line has rapidly decaying terms and 
starts with $pe^{-2^{6\log_2(r)}M^{-2}}$ which is smaller than
$pe^{-r^2}$. Thus 
\begin{equation*}
  \P\l[|a_n|\le \frac{2^k}{M} \mb{ for } n\in A_k \mb{
  for every }k>3\log_2(r)\r] = 1-o(1).
\end{equation*}

Thus the event in question has probability at least $e^{-
  cr^{2-\alp}(\log(r))^2(1+o(1))}$.



\item The event $|a_n|\le \frac{2^k}{M}$ for $n\in D_k$ 
  for $k\le \lceil \frac{M}{p}\rceil$. Following exactly the same
  steps as above we can prove that
\begin{equation*}
\P\l[|a_n|\le \frac{2^k}{M} \mb{ for } n\in D_k \r]
  \ge e^{- cr^{2-\alp}(\log(r))^2(1+o(1))}.
\end{equation*}

\item The event $\summ_{n\in B} |a_n|\frac{r^n}{\sqrt{n!}} \le 4 \frac{r^m}{\sqrt{m!}}$. 
By Cauchy-Schwarz, $\l(\summ_{n\in B} |a_n|\frac{r^n}{\sqrt{n!}}\r)^2 \le \l( \summ_{n\in B} |a_n|^2 \r) \l(\summ_{n\in B} \frac{r^{2n}}{n!} \r)$. $Y=\summ_{n\in B} |a_n|^2$ has $\Gam(|B|,1)$ distribution. Also
$\summ_{n\in B} \frac{r^{2n}}{n!} \le e^{r^2}$, since the left hand
is part of the Taylor series of $e^{r^2}$.
Therefore the event in question has probability, 
\begin{eqnarray*}
\P[\mb{event in question}] &\ge& \P\l[Y< 16 \frac{r^{2m}}{m!}e^{-r^2}\r] \\
             &\ge& \varphi\l(8 \frac{r^{2m}}{m!}e^{-r^2}\r)8 \frac{r^{2m}}{m!}e^{-r^2},
\end{eqnarray*} 
where $\varphi$ is the density of the $\Gam(|B|,1)$ distribution.  This last follows because $\varphi$ is increasing on $[0,|B|]$ and thus $\P[Y<x]\ge \varphi\l(\frac{x}{2}\r)\frac{x}{2}$. for $x<|B|$.
Continuing,

\begin{eqnarray}
\P[\mb{event in question}] &\ge&\frac{1}{(2\gam r^{\alp})!} e^{-8
  \frac{r^{2m}}{m!}e^{-r^2}}\l(8 \frac{r^{2m}}{m!}e^{-r^2}\r)^{2\gam
  r^{\alp}} \\ 
&\ge& Ce^{2\gam r^{\alp} m\log
  (r^2)-2\gam r^{2+\alp}-2\gam r^{\alp}m\log 
  (m)+2\gam r^{\alp}m +O(r^{\alp}\log(r))}. 
\label{eq:probofB} 
\end{eqnarray}
where we used Stirling's approximation.

The exponent needs simplification. Take the first and third terms in
the exponent. We have $-2\gam
mr^{\alp}\log\l(\frac{m}{r^2}\r)$. Recall that $m=r^2+\gam r^{\alp}$
and that $\alp<2$. Therefore by Taylor's expansion of $\log(1+\gam
r^{\alp-2})$ we get 
\begin{equation}\label{eq:stirl}
-2\gam mr^{\alp}\log \l(\frac{m}{r^2}\r) = \l\{  \begin{array}{ccc}
           -2\gam^2 r^{2\alp} &+ \gam^3 r^{3\alp-2} &-\frac{2}{3} \gam^4  r^{4\alp-4}+\ldots\\
           -2\gam^3 r^{3\alp-2} &+\gam^4 r^{4\alp-4} &-\frac{2}{3} \gam^5 r^{5\alp-6}+\ldots \end{array} \r.
\end{equation}
Now consider (\ref{eq:probofB}). Expand the fourth term in the 
exponential as $2\gam r^{2+\alp}+2\gam^2 r^{2\alp}$. We get the following
terms 
\begin{itemize}
\item $r^{2+\alp}(-2\gam+2\gam)=0$, from the second and fourth terms
  (first piece of the fourth term) in the exponential in (\ref{eq:probofB}).
\item $r^{2\alp}(-2\gam^2 + 2\gam^2)=0$, from the sum of the first
  term in the expansion (\ref{eq:stirl}) and the second piece of the
  fourth term in the exponential  in (\ref{eq:probofB}).  
\item $r^{3\alp-2}(\gam^3-2\gam^3)=-\gam^3r^{3\alp-2}$, from the
  expansion (\ref{eq:stirl}). 
\item Other terms such as
  $r^{\alp}\log(m),r^{\alp}\log(r),r^{\alp},r^{4\alp-4},r^{5\alp-6}$
  etc. All these are of lower order than $r^{3\alp-2}$ when $1<\alp<2$.
\end{itemize}
Hence,
\[ \P[\mb{event in quesion}] \ge e^{-\gam^3 r^{3\alp-2}(1+o(1))}. \]

\item The event $|a_n|<n-2m$ for $n\in C$. This is just an event for a sequence of i.i.d.  Complex Gaussians. It has a fixed probability $p_0$ (say). 
\item The event $|a_m| \ge 15$ also has a constant probability (not
  depending on $r$, that is).

\end{enumerate}
This completes the estimation of probabilities. Among these five
events, the third one, namely $\summ_{n\in B} |a_n|\frac{r^n}{n!} \le
4 \frac{r^m}{\sqrt{m!}}$ has the least probability (Recall that $1 <
 \alp < 2$).  

Also these events are all independent, being dependent on disjoint sets of coefficients. Thus $\P\l[n(r)\ge r^2+\gam r^{\alp} \r]\ge e^{-\gam^3r^{3\alp-2}(1+o(1))}$. 
\end{proof}

\vspace{.5 in} \noindent{\bf Acknowledgments.} I am grateful to Yuval
Peres and Mikhail Sodin for suggesting the problem of overcrowding and
the problem of moderate deviations respectively. I also thank them for
innumerable discussions which illuminated many aspects of
 two problems considered here. Earlier I had proved the lower bound  in
Theorem~\ref{thm:moderate} only for  $\alp>\frac{4}{3}$. I am very
thankful to Fedor Nazarov for showing me how  to prove (communicated
to me through Misha Sodin) the bound for all $\alp>1$.

\sc \bigskip \noindent Manjunath Krishnapur, Department of
Statistics, U.C.\ Berkeley, CA 94720,
USA. \\{\sf manju@stat.berkeley.edu}

\end{document}